\documentclass{amsart}
\usepackage{amsmath}
\usepackage{amssymb}
\usepackage{amsthm}
\usepackage[matrix,arrow,curve]{xy}
\allowdisplaybreaks
\makeatletter
\@namedef{dgo@ral}{\let\dg@VECTOR=\rightleftvector}
\def\rightleftvector(#1,#2)#3{%
   \begingroup
   \dg@XTEMP=#1\relax\multiply\dg@XTEMP\m@ne\relax
   \dg@YTEMP=#2\relax\multiply\dg@YTEMP\m@ne\relax
   \begin{picture}(0,0)%
      \thinlines
      \put(0,60){\vector(#1,#2){#3}}%
      \put(#3,-60){\vector(\dg@XTEMP,\dg@YTEMP){#3}}%
   \end{picture}%
   \endgroup}%
\@namedef{dgo@rar}{\let\dg@VECTOR=\rightrightvector}
\def\rightrightvector(#1,#2)#3{%
   \begingroup
   \dg@XTEMP=#1\relax\multiply\dg@XTEMP\m@ne\relax
   \dg@YTEMP=#2\relax\multiply\dg@YTEMP\m@ne\relax
   \begin{picture}(0,0)%
      \thinlines
      \put(0,60){\vector(#1,#2){#3}}%
      \put(0,-60){\vector(#1,#2){#3}}%
   \end{picture}%
   \endgroup}%
\makeatother
\numberwithin{equation}{section}
\newtheorem{thm}{Theorem}[section]
\newtheorem{lem}[thm]{Lemma}
\newtheorem{defn}[thm]{Definition}
\newtheorem{prop}[thm]{Proposition}

\newtheorem{rem}[thm]{Remark}

\newcommand{\G}{\Omega}
\newcommand{\widebar}[1]{\overline{#1}}

\newcommand{\ad}{\operatorname{ad}}

\newcommand{\proj}{\mbox{pr}}
\newcommand{\comp}{\circ}

\newcommand{\sq}[1]{{Sq}^{#1}}

\newcommand{\field}{\mathbb F}

\newcommand{\integral}{\mathbb Z}

\newcommand{\complex}{\mathbb C}

\newcommand{\cuplen}[1]{\operatorname{cup}(#1)}

\newcommand{\wgt}[1]{\operatorname{wgt}(#1)}

\newcommand{\Mwgt}[1]{\operatorname{{M}wgt}(#1)}

\newcommand{\Cat}[1]{\operatorname{Cat}(#1)}

\newcommand{\cat}[1]{\operatorname{cat}(#1)}

\newcommand{\sigmacat}[1]{\operatorname{\textrm{r}cat}(#1)}
\newcommand{\spin}[1]{\textrm{\bf Spin}(#1)}

\newcommand{\gtwo}[1]{\mathrm{\bf G}_{2}}
\newcommand{\ffour}[1]{\mathrm{\bf F}_{4}}
\newcommand{\esix}[1]{\mathrm{\bf E}_{6}}
\newcommand{\eseven}[1]{\mathrm{\bf E}_{7}}
\newcommand{\eeight}[1]{\mathrm{\bf E}_{8}}
\newcommand{\smashprod}{\wedge}
\newcommand{\bigsmashprod}{\bigwedge}

\newcommand{\Ext}{\operatorname{Ext}}

\newcommand{\Min}{\operatorname{Min}}
\newcommand{\Max}{\operatorname{Max}}

\makeatletter
\def\proofname{{\it Proof:}}
\def\aproofname{{\it Proof of }}
\def\pproofname{{\it. }}
\def\@proof[#1]{\aproofname{\it #1}\pproofname}
\renewenvironment{proof}{\par\noindent\hspace{18pt}\@ifnextchar[{\@proof}{{\proofname}\quad }}{{\unskip\nobreak\hfill{\it\qedsymbol}}\par\vskip 9pt}
\ifx\undefined\bysame\newcommand{\bysame}{\leavevmode\hbox to3em{\hrulefill}\,}\fi
\makeatother
\title{%
Lusternik-Schnirelmann category of $\spin{9}$%
}
\author[Iwase]{Norio Iwase%
${}^{\dagger}$%
}
\email[Iwase]{iwase@math.kyushu-u.ac.jp}
\author[Kono]{Akira Kono%
}
\email[Kono]{kono@kusm.kyoto-u.ac.jp}
\address[Iwase]{Faculty of Mathematics,
  Kyushu University,
  Fukuoka 810-8560, Japan}
\address[Kono]{Department of Mathematics,
  Faculty of Science,
  Kyoto University,
  Kyoto 607-8502, Japan}
\thanks{${}^{\dagger}$supported by the Grant-in-Aid for Scientific Research \#14654016 from Japan Society for the Promotion of Science.}
\begin{document}
%
%     ABSTRACT
%
\begin{abstract}
Let $G$ be a compact connected Lie group and $p : E \to {\Sigma}^2V$ a principal $G$-bundle with a characteristic map $\alpha : A{=}{\Sigma}V \to G$.
By combining cone decomposition arguments in \cite{IMN:spin(7),IMN:cat-lie} with computations of higher Hopf invariants introduced in \cite{Iwase:counter-ls-m}, we generalize the result in \cite{IM:cat-sp(3)}:
Let $\{F_{i}\,\vert\,0{\leq}i{\leq}m\}$ be a cone-decomposition of $G$ with a canonical structure map $\sigma_{i}$ of $\cat{F_{i}} \leq i$ for $i \leq m$.
We have $\cat{E} \leq \Max(m{+}n,m{+}2)$ for $n \geq 1$, if $\alpha$ is compressible into $F_{n} \subseteq F_{m}$ $\simeq$ $G$ and $H^{\sigma_n}_n(\alpha) = 0$, under a suitable compatibility condition.
On the other hand, calculations of \cite{HK:comm_spin} and \cite{IKT:cohomology_spin} on spinor groups yields a lower estimate for the L-S category of spinor groups by means of a new computable invariant $\Mwgt{-;{\field_2}}$ which is stronger than $\wgt{-;{\field_2}}$ introduced in \cite{Rudyak:ls-cat_mfds2} and \cite{Strom:essential-cat-wgt}.
As a result, we obtain $\cat{\spin{9}}=\Mwgt{\spin{9};{\field_2}} = 8 > 6 = \wgt{\spin{9};{\field_2}}$.
\end{abstract}
\maketitle
%
%	Document
%
\baselineskip 18pt
\section*{Introduction}
The Lusternik-Schnirelmann category $\cat{X}$, L-S category for short, is the least integer $m$ such that there is a covering of $X$ by $(m{+}1)$ open subsets each of which is contractible in $X$.
Ganea introduced a stronger notion of L-S category, $\Cat{X}$ the strong L-S category of $X$, which is equal to the cone-length by Ganea \cite{Ganea:category}, that is, the least integer $m$ such that there is a set of cofibre sequences $\{A_{i} \to X_{i-1} \hookrightarrow X_{i}\}_{1 \leq i \leq m}$ with $X_{0} = \{\ast\}$ and $X_{m} \simeq X$.
Then by Ganea \cite{Ganea:category}, we have $\cat{X} \leq \Cat{X} \leq \cat{X} {+} 1$.
Throughout this paper, we follow the notations in \cite{Iwase:counter-ls,Iwase:counter-ls-m}:
For a map $f : S^{k} \to X$, a homotopy set of higher Hopf invariants $H_m^S(f) = \{[H^{\sigma}_m(f)]\,\vert\,\text{$\sigma$ is a structure map of $\cat{X}{\leq}m$}\}$ (or its stabilisation ${\mathcal H}_m^S(f) = {\Sigma}^{\infty}_{\ast}H_m^S(f)$) is referred simply as {\it a (stabilised) higher Hopf invariant} of $f$, which plays a crucial role in this paper.

A computable lower estimate is given by the classical cup-length.
Here we give its definition in a slightly general fashion:
\begin{defn}[I. \cite{Iwase:ls-cat-survay}]
\begin{enumerate}
\item
Let $h$ be a multiplicative generalized cohomology.\par\noindent
$\begin{displaystyle}
\cuplen{X;h} = \Min \left\{ m \geq 0 \left\vert\, \forall{\{v_0,\cdots,v_m \in \tilde{h}^{\ast}(X)\}} \ v_0{\cdot}v_1{\cdot}{\cdots}{\cdot}v_m=0\right.\right\}
\end{displaystyle}$.
\item
$\begin{displaystyle}
\cuplen{X} = \Max \left\{ \cuplen{X;h} \left\vert\, \text{\begin{minipage}[c]{50mm}
$h$ is a multiplicative generalized cohomology
\end{minipage}}
\right.\right\}.
\end{displaystyle}$
\end{enumerate}
\end{defn}
Then we have $\cuplen{X;h} \leq \cuplen{X} \leq \cat{X}$ for any multiplicative generalized cohomology $h$.
When $h$ is the ordinary cohomology with a coefficient ring $R$, we denote $\cuplen{X;h}$ by $\cuplen{X;R}$.
This definition immediately implies the following.
\begin{rem}
$\begin{displaystyle}
\cuplen{X} = \Min\left\{ m \geq 0 \left\vert\, \text{$\tilde{\Delta}^{m+1} : X \to \smashprod^{m+1}X$ is stably trivial} \right.\right\}.
\end{displaystyle}$
\end{rem}

Let $\{p^{{\G}X}_k \,{:}\, E^{k}({\G}{X}) {\to} P^{k-1}({\G}X) \,;\, k {\geq} 1\}$ be the $A_{\infty}$-structure of ${\G}X$ in the sense of Stasheff \cite{Stasheff:higher-associativity} (see also Iwase-Mimura \cite{IM:higher-associativity} for some more properties).
The relation between an $A_{\infty}$-structure and a L-S category gives the key observation in \cite{Iwase:counter-ls,Iwase:counter-ls-m,Iwase:ls-cat-bundle} to producing counter-examples to the Ganea conjecture on L-S category.
On the other hand, Rudyak \cite{Rudyak:ls-cat_mfds2} and Strom \cite{Strom:essential-cat-wgt} introduced a homotopy theoretical version of Fadell-Husseini's category weight (see \cite{FH:category_weight}), which can be described as follows, for an element $u \in h^{\ast}(X)$ and a generalized cohomology $h$:
\begin{align*}&
\wgt{u;h} = \Min\left\{m \geq 0%
\left\vert\,%
(e^{X}_{m})^{\ast}(u) \not= 0 \ \text{in} \ h^{\ast}(P^m(\G{X}))
\right.\right\},
\\&
\wgt{X;h} = \Min\left\{m \geq 0%
\left\vert\,%
\begin{array}[c]{l}
\text{\begin{minipage}[c]{50mm}
$(e^{X}_{m})^{\ast} : h^{\ast}(X) \to h^{\ast}(P^m(\G{X}))$ is a monomorphism
\end{minipage}}
\end{array}
\right.\right\},
\end{align*}
where $e^{X}_{m}$ denotes the map $P^{m}(\G{X}) \hookrightarrow P^{\infty}(\G{X}) \simeq X$.
Then we easily see 
\begin{equation}\label{eq:weight}
\wgt{X;h} = \Max\{\wgt{u;h} \,\vert\, u \in \tilde{h}^{\ast}(X)\}.
\end{equation}
We remark that $\wgt{u;h} = s$ if and only if $u$ represents a non-zero class in $E^{s,\ast}_{\infty}$ of bar spectral sequence $\{(E^{\ast,\ast}_r,d^{\ast,\ast}_r)\,\vert\,r\geq1\}$ converging to $h^{\ast}(X)$ with $E_{2}^{**} \cong \Ext^{**}_{h_{(\G{X})}}(h_{*},h_{*})$.
When $h$ is the ordinary cohomology with a coefficient ring $R$, we denote $\wgt{X;h}$ by $\wgt{X;R}$.
In this paper, we introduce new computable invariants in terms of (unstable) cohomology operations:
\begin{defn}
Let $h$ be a generalized cohomology:
A homomorphism $\phi : h^*(X) \to h^*(Y)$ is called a $h$-morphism if it preserves the actions of all (unstable) cohomology operations on $h^*$.
\end{defn}
\begin{defn}[I. \cite{Iwase:ls-cat-survay}]
Let $h$ be a generalized cohomology and $X$ a space.
A module weight $\Mwgt{X;h}$ of $X$ with respect to $h$ is defined as follows:
$$
\Mwgt{X;h} = \Min\left\{ m {\geq} 0 \left\vert\, \text{\begin{minipage}[c]{75mm}There is an $h$-morphism $\phi$ : $h^{\ast}(P^m(\G{X}))$ $\to$ $h^{\ast}(X)$, which is a left homotopy inverse of $(e^{X}_{m})^{\ast}$. 
\end{minipage}}
\right.\right\}
$$
\end{defn}
When $h$ is the ordinary cohomology with coefficients in a ring $R$, we denote $\Mwgt{X;h}$ by $\Mwgt{X;R}$.
These invariants satisfy the following inequalities:
\begin{align*}
\cuplen{X;R} \leq \wgt{X;R} \leq &\Mwgt{X;R} \leq \cat{X}.
\end{align*}
Similar to the above definition of $\cuplen{X}$, we define the following invariants:
\begin{defn}[I. \cite{Iwase:ls-cat-survay}]
\begin{enumerate}
\item
$\begin{displaystyle}
\wgt{X} = \Max\left\{\wgt{X;h}\left\vert\,%
\text{\begin{minipage}[c]{28mm}
$h$ is a generalized cohomology
\end{minipage}}
\right.\right\}
\end{displaystyle}$
\item
$\begin{displaystyle}
\Mwgt{X} = \Max\left\{\Mwgt{X;h}\left\vert\,%
\text{\begin{minipage}[c]{45mm}
$h$ is a generalized cohomology
\end{minipage}}
\right.\right\}
\end{displaystyle}$
\end{enumerate}
\end{defn}
\begin{rem}
Let $\sigmacat{-}$ be Rudyak's stable L-S category, which is denoted as $r(-)$ in \cite{Rudyak:ls-cat_mfds2}.
Then we have 
$\cuplen{X} \leq \wgt{X} = \sigmacat{X} \leq \Mwgt{X} \leq \cat{X}.$
\end{rem}

\section{Main results}\label{sect:1}

From now on, we work in the category of connected CW-complexes and continuous maps.
Let us denote by $Z^{(k)}$ the $k$-skeleton of a CW complex $Z$.
To give an upper-bound for L-S category of the total space of a fibre bundle $F \hookrightarrow E \to B$, we need a refinement of results of Varadarajan \cite{Varadarajan:fib-cat} and Hardie  \cite{Hardie:fib-cat}, and corresponding result for strong category of Ganea \cite{Ganea:category}:
\begin{thm}[\cite{Varadarajan:fib-cat,Hardie:fib-cat,Ganea:category}]\label{thm:ganea}
\begin{enumerate}
\item~
$\cat{E}{+}1 \leq (\cat{F}{+}1){\cdot}(\cat{B}{+}1)$ 
\item~
$\Cat{E}{+}1 \leq (\Cat{F}{+}1){\cdot}(\Cat{B}{+}1)$.
\end{enumerate}
\end{thm}
In \cite{IMN:cat-lie}, Iwase-Mimura-Nishimoto gave such a refinement in the case when the base space $B$ is non-simply connected.
But in this paper, we give another refinement in the case when the fibre bundle is a principal bundle over a double suspension space:
Let $G$ be a compact Lie group with a cone-decomposition of length $m$: 
$$
\text{($m$ cofibre sequences)}\quad K_{i} \overset{\rho_{i}}\to F_{i{-}1} \hookrightarrow F_{i},\quad i \geq 1,
$$
with $F_{0} = \{\ast\}$ and $F_{i} = F_{m} \simeq G$, $i \geq m$.
Then we obtain $\sigma_k : F_{k} \to P^k\G{F_{k}}$ for all $k \leq m$ as a right homotopy inverse of $e_k : P^k\G{F_{k}} \to P^\infty\G{F_{k}} \simeq F_{k}$ by induction on $k \geq 1$.
Thus we have the following commutative diagram:
$$
\xymatrix{
\ \{\ast\} \ 
\ar[d]_{\sigma_0}
\ar@{^{(}->}[r]
&
\ F_1 \ 
\ar[d]_{\sigma_1}
\ar@{^{(}->}[r]
&
\ \cdots \ 
\ar@{^{(}->}[r]
&
\ F_m \ 
\ar[d]_{\sigma_m}
\\
\ \{\ast\} \ 
\ar@{^{(}->}[r]
\ar[d]_{e_0}
&
\ P^1\G{F_1} \ 
\ar@{^{(}->}[r]
\ar[d]_{e_1}
&
\ \cdots \ 
\ar@{^{(}->}[r]
&
\ P^m\G{F_m} \ 
\ar[d]_{e_m}
\\
\ \{\ast\} \ 
\ar@{^{(}->}[r]
&
\ F_1 \ 
\ar@{^{(}->}[r]
&
\ \cdots \ 
\ar@{^{(}->}[r]
&
\ F_m, \ 
}
$$
where $e_k{\comp}\sigma_k \sim 1_{F_k}$ for all $k \leq m$.
\begin{thm}\label{thm:main}
Let $G \hookrightarrow E \to {\Sigma}^2V$ be a principal bundle with characteristic map $\alpha : A{=}{\Sigma}V \to G$.
Then we have $\cat{E} \leq \Max(m{+}n,m{+}2)$ for $n \geq 1$, if 
\begin{enumerate}
\item\label{ass:compatible-decomposition3}
$\alpha$ is compressible into $F_{n} \subseteq F_{m}$ $\simeq$ $G$, 
\item\label{ass:compatible-decomposition2}
$H^{\sigma_n}_n(\alpha) = 0$ and
\item\label{ass:compatible-decomposition1}
the restriction of the multiplication $\mu : G{\times}G \to G$ to $F_{j}{\times}F_{n} \subseteq F_{m}{\times}F_{m}$ $\simeq$ $G{\times}G$ is compressible into $F_{j{+}n} \subseteq F_{m}$ $\simeq$ $G$, $j {\geq} 0$ as $\mu_{j,n} : F_{j}{\times}F_{n} \to F_{j{+}n}$ such that $\mu_{j,n}\vert_{F_{j-1}{\times}F_{n}} = \mu_{j-1,n}$.
\end{enumerate}
\end{thm}
\begin{rem}\label{rem:cellular}
If we choose $n=m{+}1$, then the assumptions (\ref{ass:compatible-decomposition3}) through (\ref{ass:compatible-decomposition1}) above are automatically satisfied.
Thus we always have $\Cat{E} \leq 2\Cat{G}{+}1$ which is a special case of Ganea's theorem (see Theorem \ref{thm:ganea} (2)).
\end{rem}

For $\spin{9}$, we first observe the following result.
\begin{thm}\label{thm:main1}
$\Mwgt{\spin{9};{\field_2}} \geq 8 > 6 = \wgt{\spin{9};{\field_2}}$.
\end{thm}
These results imply the following result.
\begin{thm}\label{thm:main2}
$\begin{textstyle}
\cat{\spin{9}}=\Mwgt{\spin{9};{\field_2}}=8.
\end{textstyle}$
\end{thm}
\section{Proof of Theorem \ref{thm:main}}\label{sect:2}
From now on, we work in the homotopy category, and so we do not distinguish $G$ and $F_m$.
Let $G \hookrightarrow E \to {\Sigma}^2V$ be a principal bundle with characteristic map $\alpha : A{=}{\Sigma}V \to G$.
The assumptions (\ref{ass:compatible-decomposition3}) and (\ref{ass:compatible-decomposition1}) in Theorem \ref{thm:main} allows us to construct a filtration $\{E_{k}\}_{0 \leq k \leq n{+}m}$ of $E$:
By using the James-Whitehead decomposition (see Theorem 1.15 of Whitehead \cite{Whitehead:elements}), we have
\begin{align*}&\begin{textstyle}
E = G \cup_{\psi} G{\times}CA,\quad 
\psi=\mu{\circ}(1_G{\times}\alpha) : G{\times}A \to G.
\end{textstyle}\end{align*}
Firstly, we define two filtrations of $E$ as follows:
\begin{align*}&
E_{k} = 
\left\{\begin{array}[c]{ll}
F_{k},
&\quad k \leq n,
\\[1mm]
F_{k} \cup_{\psi_{k-n-1}} F_{k-n-1} {\times} C{A},
&\quad n < k \leq m{+}n,
\end{array}\right.
\\&
E'_{k} = 
\left\{\begin{array}[c]{ll}
F_{k},
&\quad k < n,
\\[1mm]
F_{k} \cup_{\psi_{k-n}} F_{k-n} {\times} C{A},
&\quad n \leq k \leq m{+}n,
\end{array}\right.
\end{align*}
where $\psi_{j} = \mu_{j,n} {\circ} (\alpha{\times}1) : F_{j}{\times}A \to F_{j+n}$ and $E=E'_{m+n}$ which immediately imply that $\cat{E} = \cat{E'_{m+n}}$.

By using the assumption (\ref{ass:compatible-decomposition1}), we obtain the following cofibre sequences, similarly to the arguments given in \cite{IMN:cat-lie}:
\begin{align*}&
K_{k+1} \to E_{k} \hookrightarrow E_{k+1}, \ \text{for $0 \leq k < n$,}
\\[1mm]&
K_{k+1} \vee (K_{k-n}{\ast}A) \to E_{k} \hookrightarrow E_{k+1}, \ \text{for $n \leq k < m{+}n$,}
\\[1mm]&
K_{k-n}{\ast}A \to E_{k} \hookrightarrow E'_{k},
\end{align*}
Similarly to the arguments given in \cite{IMN:spin(7),IMN:cat-lie}, we obtain 
\begin{equation}\label{eq:imn}
\cat{E_{k}} \leq k \  \text{and} \ \cat{E'_{k}} \leq k{+}1\quad \text{for any $k \geq n$},
\end{equation}
by induction on $k$.
The following lemma can be deduced in a similar but easier manner to the main theorem of \cite{Iwase:ls-cat-bundle}, using $H^{\sigma_n}_n(\alpha) = 0$, the assumption (\ref{ass:compatible-decomposition2}):
\begin{lem}\label{lem:key}
$\cat{E'_{j+n}} \leq j{+}n$ for all $j \geq 0$ and $n \geq 2$.
\end{lem}
Lemma \ref{lem:key} and (\ref{eq:imn}) imply $\cat{E} = \cat{E'_{m{+}n}} \leq \Max(m{+}n,m{+}2)$, and hence we are left to show Lemma \ref{lem:key}.
\begin{proof}[Lemma \ref{lem:key}]
We define a map $\hat{\psi}_{j}$ as follows:
$$
\hat{\psi}_{j} = \sigma_{j+n}{\comp}\mu_{j,n}{\comp}(e_{j}{\times}\alpha) : P^j\G{F_j}{\times}A \to P^{j+n}\G{F_{j+n}}.
$$
Then we have $\hat{\psi}_{j}{\comp}(\sigma_j{\times}1) = \sigma_{j+n}{\comp}\mu_{j,n}{\comp}(e_{j}{\times}\alpha){\comp}(\sigma_j{\times}1) \sim \sigma_{j+n}{\comp}\mu_{j,n}{\comp}(1{\times}\alpha) = \sigma_{j+n}{\comp}\psi_{j}$ and $e_{j+n}{\comp}\hat{\psi}_{j} = e_{j+n}{\comp}\sigma_{j+n}{\comp}\mu_{j,n}{\comp}(e_{j}{\times}\alpha) \sim \mu_{j,n}{\comp}(e_{j}{\times}\alpha) = \psi_j{\comp}(e_{j}{\times}1)$.
Thus the following diagram is commutative up to homotopy:
\begin{equation}\label{eq:pushpull_3}
\begin{xy}
\xymatrix{
\ F_{j} \ 
	\ar[d]_{\sigma_{j}}
&
\ F_{j}{\times}A \ 
	\ar[l]_{\proj_1}
	\ar[r]^{\psi_j}
	\ar[d]_{\sigma_{j}{\times}1}
&
\ F_{j+n} \ 
	\ar[d]^{\sigma_{j+n}}
\\
\ P^j\G{F_j} \ 
	\ar[d]_{e_j}
&
\ P^j\G{F_j}{\times}A \ 
	\ar[l]_{\proj_1}
	\ar[r]^{\hat\psi_j}
	\ar[d]_{e_{j}{\times}1}
&
\ P^{j+n}\G{F_{j+n}} \ 
	\ar[d]^{e_{j+n}}
\\
\ F_{j} \ 
&
\ F_{j}{\times}A \ 
	\ar[l]_{\proj_1}
	\ar[r]^{\psi_j}
&
\ F_{j+n} \ 
}
\end{xy}
\end{equation}
Therefore, the space $E'_{j+n} = F_{j+n} \cup_{\psi_j}F_j{\times}CA$ is dominated by $P^{j+n}\G{F_{j+n}} \cup_{\hat\psi_j}P^j\G{F_j}{\times}CA$.
Since $\alpha$ satisfies $H^{\sigma_n}_n(\alpha) = 0$, we have the following commutative diagram up to homotopy:
$$
\begin{xy}
\xymatrix{
\ A \ 
\ar[r]^{\alpha}
\ar[d]_{\Sigma\ad{\alpha}}
&
\ F_m \ 
\ar[d]^{\sigma_n}
\\
\ \Sigma\G{F_n} \ 
\ar@{^{(}->}[r]
&
\ P^n\G{F_n}, \ 
}
\end{xy}
$$
where $\ad{\alpha} : V \to \G{F_n}$ is the adjoint map of $\alpha : A{=}{\Sigma}V \to F_n$.
Thus $\sigma_n{\comp}\alpha$ is compressible into $\Sigma\G{F_n}$.
Since $\Cat{P^i\G{F_j}{\times}\Sigma\G{F_n}} \leq i{+}1$ for $i \leq j$, we have $\hat\psi_j \vert_{P^i\G{F_j}{\times}\Sigma\G{F_n}}$ can be compressible into $P^{i+1}\G{F_{j+n}}$ for $i \leq j$.
Thus we have the following cone decomposition of $P^{j+n}\G{F_{j+n}} \cup_{\hat\psi_j}P^j\G{F_j}{\times}CA$:
\begin{align*}&
\G{F_{j+n}} \to \{\ast\} \hookrightarrow P^{1}\G{F_{j+n}},
\\&
E^{2}\G{F_{j+n}} {\vee} A \to P^{1}\G{F_{j+n}} \hookrightarrow P^{2}\G{F_{j+n}} {\cup} CA,
\\&
\quad\vdots
\\&
E^{i}\G{F_{j+n}} {\vee} E^{i-2}\G{F_{j}}{\ast}A \begin{array}[t]{l}\to P^{i-1}\G{F_{j+n}} {\cup} P^{i-3}\G{F_{j}}{\times}CA \\[0mm]\qquad\qquad\qquad\qquad\hookrightarrow P^{i}\G{F_{j+n}} {\cup} P^{i-2}\G{F_{j}}{\times}CA,\end{array}
\\&
\quad\vdots
\\&
E^{j+2}\G{F_{j+n}} {\vee} E^{j}\G{F_{j}}{\ast}A \begin{array}[t]{l}\to P^{j+1}\G{F_{j+n}} {\cup} P^{j-1}\G{F_{j}}{\times}CA \\[0mm]\qquad\qquad\qquad\qquad\hookrightarrow P^{j+2}\G{F_{j+n}} {\cup} P^{j}\G{F_{j}}{\times}CA,\end{array}
\\&
\quad\vdots
\\&
E^{j+n}\G{F_{j+n}} \begin{array}[t]{l}\to P^{j+n-1}\G{F_{j+n}} \cup_{\hat\psi_j}P^{j}\G{F_j}{\times}CA \hookrightarrow P^{j+n}\G{F_{j+n}} \cup_{\hat\psi_j}P^j\G{F_j}{\times}CA.\end{array}
\end{align*}
This implies $\Cat{P^{j+n}\G{F_{j+n}} \cup_{\hat\psi_j}P^j\G{F_j}{\times}CA} \leq j{+}n$ for all $j \geq 0$ and $n \geq 2$, and hence $\cat{E'_{j+n}} = \cat{F_{j+n} \cup F_j{\times}CA} \leq j{+}n$ for all $j \geq 0$ and $n \geq 2$.
\end{proof}
This completes the proof of Theorem \ref{thm:main}.
\section{Bar spectral sequence}\label{sect:3}
To calculate our module weight $\Mwgt{X;{\field_2}}$ together with $\wgt{X;{\field_2}}$, we need to know the module structure of $H^*(P^{m}(\spin{9});\field_2)$ over the Steenrod algebra modulo $2$.
By Ishitoya-Kono-Toda \cite{IKT:cohomology_spin}, Hamanaka-Kono \cite{HK:comm_spin} and Kono-Kozima \cite{KK:spinor},  the followings are known:
\begin{align*}&
H^{\ast}(\spin{9};{\field_2}) = {\field_2}[x_3]/(x^4_3){\otimes}{\wedge}_{\field_2}(x_5,x_7,x_{15}), \ 
\\&\qquad
\sq{2}x_3=x_5, \ \sq{1}x_5=x_6, \ \ x_i \in H^{i}(\spin{9};{\field_2}),
\\[2mm]&
H_{\ast}(\G\spin{9};{\field_2}) = \wedge_{\field_2}(u_2){\otimes}{\field_2}[u_4,u_6,u_{10},u_{14}], \ 
\\&\qquad
u_4\sq{2}=u_2, \ u_{10}\sq{2}=u^{2}_{4}, \ u_{14}\sq{4}=u_{10}, \ \ u_{2i} \in H_{2i}(\G\spin{9};{\field_2}),
\end{align*}
where we denote by $\wedge_R(a_{i_1},\cdots,a_{i_t})$ the exterior algebra on $a_{i_1},\cdots,a_{i_t}$ over $R$.
We remark that the cohomology suspension of $x_{2i+1}$ is non-trivially given by $u_{2i}$ for $i=1$, $2$, $3$ and $7$.
To determine $H^*(P^{m}(\spin{9});\field_2)$, we have to study the bar spectral sequence $(E^{*,*}_r,d^{*,*}_r)$ converging to $H^*(\spin{9};\field_2)$:
\begin{align*}&
E^{s,t}_1 \cong \tilde{H}^{s{+}t}(P^{s}(\G\spin{9}),P^{s{-}1}(\G\spin{9});\field_2) 
\cong \tilde{H}^{t}(\bigwedge^s\G\spin{9};\field_2),
\\[2mm]&
D^{s,t}_1 \cong \tilde{H}^{s{+}t}(P^{s}(\G\spin{9});\field_2),
\\[2mm]&
E^{*,*}_2 \cong \Ext^{*,*}_{H_{\ast}(\G\spin{9};{\field_2})}({\field_2},{\field_2}) 
\cong \field_2[x_{1,2}]{\otimes}{\wedge_{\field_2}}(x_{1,4},x_{1,6},x_{1,10},x_{1,14}), 
\\[2mm]&
E^{*,*}_{\infty} \cong \tilde{H}^{\ast}(\spin{9};\field_2) \cong \field_2[x_{1,2}]/(x^4_{1,2}){\otimes}{\wedge_{\field_2}}(x_{1,4},x_{1,6},x_{1,14}),
\end{align*}
where $x_{1,2}$, $x_{1,4}$, $x_{1,6}$ and $x_{1,14}$ are permanent cycles by \cite{KK:spinor}.
Therefore, there is only one differential $d_a(x_{1,10})$ ($a \geq 2$) which is possibly non-trivial, and we have $E^{*,*}_a \cong E^{*,*}_2$ and $E^{*,*}_{a{+}1} \cong E^{*,*}_{\infty}$.
Since $x_3$ is of height 4 in $H^{\ast}(\spin{9};\field_2)$, we have $d_a(x_{1,10})=x^4_{1,2}$, and hence $a=3$.
Thus we have the following:
\begin{align*}&
d_r=0 \ \text{if $r\not=3$}, \ d_3(x_{1,i})=0 \ \text{if $i\not=10$}, \ d_3(x_{1,10})=x^{4}_{1,2},
\\[2mm]&
E^{*,*}_2 \cong E^{*,*}_3 \cong \field_2[x_{1,2}]{\otimes}{\wedge_{\field_2}}(x_{1,4},x_{1,6},x_{1,10},x_{1,14}), 
\\[2mm]&
E^{*,*}_4 \cong E^{*,*}_{\infty} \cong \field_2[x_{1,2}]/(x^4_{1,2}){\otimes}{\wedge_{\field_2}}(x_{1,4},x_{1,6},x_{1,14}).
\end{align*}
By truncating the above computations with the same differential $d_{r}$ to the spectral sequence for $P^{m}(\spin{9})$ of Stasheff's type (similar to the computation in \cite{Iwase:k-ring}),
 we are lead to the following proposition, and we leave the details to the reader.
\begin{prop}\label{prop:ring-structure}
Let $A = \field_2[x_{3}]/(x^4_{3}){\otimes}{\wedge_{\field_2}}(x_{5},x_{7},x_{15})$.
Then for $m \geq 0$, we have 
\begin{equation*}
H^*(P^m(\spin{9});\field_2) \cong 
\begin{cases}
A^{[0]} \cong \field_2, \ & \text{if $m=0$},\\
A^{[m]} \oplus x_{11}{\cdot}A^{[m-1]} \oplus S_m,\ & \text{if $3 \geq m \geq 1$},\\
A^{[m]} \oplus x_{11}{\cdot}(A^{[m{-}1]}/A^{[m{-}4]}) \oplus S_{m},\ & \text{if $m \geq 4$}
\end{cases}
\end{equation*}
as modules, where $A^{[m]}$ ($m \geq 0$) denotes the quotient module $A/D^{m{+}1}(A)$ of $A$ by the submodule $D^{m{+}1}(A) \subseteq A$ generated by all the products of $m{+}1$ elements in positive dimensions, $x_{11}{\cdot}(A^{[m{-}1]}/A^{[m{-}4]})$ ($m\geq 4$) denotes a submodule corresponding to a submodule in $\field_2[x_{3}]/(x^4_{3}){\otimes}{\wedge_{\field_2}}(x_{5},x_{7},x_{11},x_{15})$ and $S_{m}$ satisfies $S_{m}{\cdot}\tilde{H}^*(P^{m}(\spin{9});\field_2)=0$ and $S_{m}\vert_{P^{m-1}(\spin{9})}=0$.
\end{prop}
Some more comments might be required to the second direct summand of the above expressions of $H^*(P^{m}(\spin{9});\field_2)$, $m \geq 4$.
The multiplication with $x_{11}$ is a fancy notation to describe the module basis and not a usual product.
However, we may regard it is a {\it partial product} in the sense introduced in the next section.
\section{Partial products}\label{sect:3h}

Since a diagonal map $\begin{displaystyle}\Delta^{\G{X}}_n = \G{(\Delta^{X}_n)} : \G{X} \to \prod^n\G{X}=\G{(\prod^nX)}\end{displaystyle}$ is a loop map, it induces a map of projective spaces:
$$
P^m(\Delta^{\G{X}}_n) : P^m(\G{X}) \to P^m(\G{(\prod^nX)}),
$$
such that $e^{\prod^nX}_m{\circ}P^m(\Delta^{\G{X}}_n) \sim \Delta^{X}_n{\circ}e^{X}_m$.
As is seen in the proof of Theorem 1.1 in \cite{Iwase:counter-ls}, there is a natural map \begin{align*}
\varphi^{X}_m : P^m(\G{(\prod^nX)}) \to 
&\bigcup_{i_1+\cdots+i_n=m}P^{i_1}(\G{X}){\times}\cdots{\times}P^{i_n}(\G{X})
\\&\qquad\subset
P^{m}(\G{X}){\times}\cdots{\times}P^{m}(\G{X})
\end{align*}
such that $\begin{displaystyle}(e^{X}_{m}{\times}\cdots{\times}e^{X}_{m}){\circ}\varphi^{X}_m = e^{\prod^nX}_m\end{displaystyle}$.
Let $\Delta^{X,m}_n = \varphi^{X}_m{\circ}P^m(\Delta^{\G{X}}_n)$, which we call the $n$-th partial diagonal of $X$ of height $m$, or simply a {\it partial diagonal} 
\begin{align*}
\Delta^{X,m}_n : P^m(\G{X}) \to &\bigcup_{i_1+\cdots+i_n=m}P^{i_1}(\G{X}){\times}\cdots{\times}P^{i_n}(\G{X})
\\&\qquad\subset
P^{m}(\G{X}){\times}\cdots{\times}P^{m}(\G{X})
\end{align*}
such that 
$\begin{displaystyle}(e^{X}_{m}{\times}\cdots{\times}e^{X}_{m}){\circ}\Delta^{X,m}_n \sim \Delta^{X}_n{\circ}e^{X}_m\end{displaystyle}$.
This partial diagonal also yields the reduced version 
\begin{align*}
\widebar{\Delta}^{X,m}_n : P^m(\G{X}) \to &\bigcup_{i_1+\cdots+i_n=m}P^{i_1}(\G{X}){\smashprod}\cdots{\smashprod}P^{i_n}(\G{X})
\\&\qquad\subset
P^{m-n+1}(\G{X}){\smashprod}\cdots{\smashprod}P^{m-n+1}(\G{X})
\end{align*}
such that 
$\begin{displaystyle}(e^{X}_{m-n+1}{\smashprod}\cdots{\smashprod}e^{X}_{m-n+1}){\circ}\widebar{\Delta}^{X,m}_n \sim \widebar{\Delta}^{X}_n{\circ}e^{X}_m\end{displaystyle}$, where $\widebar{\Delta}^{X}_n : X \to \bigsmashprod^n X$ denotes the reduced diagonal.
Let us call $\widebar{\Delta}^{X,m}_n$ the $n$-th reduced partial diagonal of $X$ of height $m$, or simply a {\it reduced partial diagonal}.

As is well-known, the product of a multiplicative generalized cohomology $h$ is given by (reduced) diagonal, i.e.,
$$
v_1{\cdot}\cdots{\cdot}v_n = (\widebar{\Delta}^{X}_n)^{\ast}(v_1{\otimes}\cdots{\otimes}v_n) \in \bar{h}^*(X),\quad\text{for any $v_1,\cdots,v_n \in \bar{h}^*(X)$,}
$$
where $\bar{h}$ denotes the reduced cohomology associated with $h$.
So it is natural to define a `partial' product as the following way:
\begin{defn}~
For any elements $v_{1},\cdots,v_{n} \in \bar{H}^*({\Sigma}\G{X};\field_2)$ which are restrictions of elements in $\bar{H}^*(P^{m-n+1}(\G{X});\field_2)$, we define a partial product $v_1{\cdot}\cdots{\cdot}v_{n} = (\widebar{\Delta}^{\G{X},m}_n)^{\ast}(v_1{\otimes}\cdots{\otimes}v_{n})$ in $\bar{H}^*(P^{m}(\G{X});\field_2)$.
\end{defn}
\begin{rem}
Since $x_{11}$ can be extended to  an element in $\bar{H}^*(P^{3}(\G{\spin{9}});\field_2)$,  we have partial products $x_{11}{\cdot}v_1{\cdot}\cdots{\cdot}v_{n-1} = (\widebar{\Delta}^{\spin{9},m}_{n-1})^{\ast}(x_{11}{\otimes}v_1{\otimes}\cdots{\otimes}v_{n-1})$ for any elements $v_{1},\cdots,v_{n-1} \in \bar{H}^*(P^{3}(\G{\spin{9}});\field_2)$, $m{-}2 \leq n \leq m$.
In the direct sum decomposition of $H^*(P^{m}(\G{\spin{9}});\field_2)$ given in Proposition \ref{prop:ring-structure}, the direct summand $x_{11}{\cdot}(A^{[m{-}1]}/A^{[m{-}4]})$ is generated by such partial products.
\end{rem}
\section{Proof of Theorem \ref{thm:main1}}\label{sect:4}
We know $x^3_{3}x_{5}x_{7}x_{15}$ and $x_{11}{\cdot}x^3_{3}x_{5}x_{7}$ exist non-trivially in $H^{\ast}(P^{8}(\G\spin{9});\field_2)$ but $x_{11}{\cdot}x^3_{3}x_{5}x_{7}$ does not exist in $H^{\ast}(P^{9}(\G\spin{9});\field_2)$ by Proposition \ref{prop:ring-structure}.
To observe what happens on the element $x_{11}{\cdot}x^3_3x_5x_{7}$ in $H^{\ast}(P^{8}(\G\spin{9});\field_2)$, we must recall the bar spectral sequence $(E^{*,*}_r,d^{*,*}_r)$:
\begin{align*}&
[(p^{\G\spin{9}}_{9})^{\ast}(x_{11}{\cdot}x^3_3x_5x_{7})]
=d_{3}(x^3_{1,2}x_{1,4}x_{1,6}x_{1,10})
= \pm x^{7}_{1,2}x_{1,4}x_{1,6} \not= 0 \ \text{in} \ E^{*,*}_3,
\end{align*}
where we denote by $[\beta]$ the corresponding class in $E^{*,*}_3$ to an element $\beta \in E^{*,*}_1$.
Thus $(p^{\G\spin{9}}_{9})^{\ast}(x^3_3x_5x_{7}x_{11}) \not=0 \ \text{in} \ E^{9,*}_1=\tilde{H}^{\ast}(\bigwedge^{9}\G\spin{9};\field_2),
$ and hence $x_{11}{\cdot}x^3_3x_5x_{7}$ does not exist in $\tilde{H}^{\ast}(P^{9}(\G\spin{9});\field_2)$, but does in $\tilde{H}^{\ast}(P^{8}(\G\spin{9});\field_2)$.

By \cite{KK:spinor}, we know $\sq{4}(x_{11})=x_{15}$ in $H^*(P^{1}(\spin{9});\field_2)$, and hence $\sq{4}(x_{11})=x_{15}$ modulo $S_3$ in $H^*(P^{3}(\spin{9});\field_2)$ for dimensional reasons.
Thus we have 
\begin{align}\label{eq:square4-1}&
\sq{4}(x_{11}{\cdot}x^3_3x_5x_{7}) = x^3_3x_5x_{7}x_{15},\quad\text{in $H^{\ast}(P^{7}(\G\spin{9});\field_2)$.}
\\\label{eq:square4-2}&
\sq{4}(x_{11}{\cdot}x^3_3x_5x_{7}) = x^3_3x_5x_{7}x_{15} + w,\quad\text{$w \in S_{8}$ in $H^{\ast}(P^{8}(\G\spin{9});\field_2)$.}
\end{align}
The equation (\ref{eq:square4-1}) implies that any left inverse epimorphism of $(e^{\spin{9}}_{7})^{*}$ 
$$
\phi : H^*(P^7(\G\spin{9});{\field_2}) \to H^*(\spin{9};{\field_2})
$$
does not preserve the action of the modulo $2$ Steenrod operations: If such a epimorphism $\phi$ did preserve  the action of the modulo $2$ Steenrod operations, the element $\phi(x^3_3x_5x_{7}x_{15})=x^3_3x_5x_{7}x_{15}$ in $H^*(\spin{9};{\field_2})$ should lie in the image of $\sq{4}$, since $x^3_3x_5x_{7}x_{15}$ lies in the image of $\sq{4}$ in $H^*(P^7(\G\spin{9});{\field_2})$ by (\ref{eq:square4-1}).
It contradicts to the fact that $H^{32}(\spin{9};{\field_2})=0$.
Thus we have
$\Mwgt{\spin{9};{\field_2}} \geq 8$.

On the other hand, we can easily see that each generator of $H^{*}(\spin{9};\field_2) \cong \field_2[x_{3}]/(x^4_{3}){\otimes}{\wedge_{\field_2}}(x_{5},x_{7},x_{15})$ has category weight $1$, and hence by  (\ref{eq:weight}), we have $\wgt{\spin{9};{\field_2}} = 6$.
This completes the proof of Theorem \ref{thm:main1}.
\section{Proof of Theorem \ref{thm:main2}}\label{sect:5}
By \cite{IMN:spin(7)}, we can easily see that $\spin{7}$ admits a cone decomposition which satisfies the condition \ref{ass:compatible-decomposition1} in Theorem \ref{thm:main}.
Since $x_{15} \in H^{15}(\spin{7});{\field_2})$ is the modulo $2$ reduction of a generator of $H^{15}(\spin{7};{\integral}) \cong {\integral}$, the image of the attaching map $\alpha$ of the $15$-cell corresponding to $x_{15}$ must lie in $\spin{7}^{(13)}$ the $13$-skeleton of $\spin{7}$, where $\spin{7}^{(13)}$ is contained in $F_3(\spin{7})$.
To observe that the attaching map $\alpha$ satisfies the condition of Theorem \ref{thm:main} with $n=3$, we need to show that $H^{\sigma_3}_3(\alpha) = 0$.
Then we obtain $\cat{\spin{9}} \leq \Cat{\spin{7}}{+}n=5{+}3=8$ by Theorem \ref{thm:main}, while we know $\Mwgt{\spin{9};{\field_2}} \geq 8$ by Theorem \ref{thm:main1}, and hence  
$$
\cat{\spin{9}} = \Mwgt{\spin{9};{\field_2}} = 8.
$$
Let $\sigma_3$ : $F_3(\spin{7})$ $\to$ $P^3(\G{F_3(\spin{7})})$ be the canonical structure map of $\cat{F_3(\spin{7})} = 3$.
Then we are left to show that $H^{\sigma_3}_3(\alpha)=0$.
By definition, 
$$
H^{\sigma_3}_3(\alpha) : S^{14} \to E^{4}(\G{F_3(\spin{7})}),
$$
where $F_3(\spin{7})=\gtwo{}^{(11)} \cup_{{\Sigma}{\complex}P^{2}} {\Sigma}{\complex}P^{3} \cup \text{(higher cells $\geq$ $8$)}$.
Since $\G{\gtwo{}^{(11)}}$ has the homotopy type of 
${\complex}P^2 \cup \text{(higher cells $\geq$ $6$)}$, we know $\G{F_3(\spin{7})}$ has the homotopy type of ${\complex}P^{3} \cup \text{(higher cells $\geq$ $6$)}$.
Thus we observe that $E^{4}(\G{F_3(\spin{7})})$ has the homotopy type of 
$$
{\Sigma^3}{\complex}P^3{\wedge}S^2{\wedge}S^2{\wedge}S^2
 \cup 
{\Sigma^3}{\complex}P^2{\wedge}{\complex}P^2{\wedge}{\complex}P^2{\wedge}{\complex}P^2
 \cup 
 \text{(higher cells $\geq$ $15$)}.
$$
It is well-known that ${\Sigma}{\complex}P^3 = {\Sigma}{\complex}P^2 \cup_{\omega_{3}} e^{7}$, $\omega_{3} : S^{6} \to S^{3} \subset {\Sigma}{\complex}P^3$, and hence we have ${\Sigma^3}{\complex}P^3{\wedge}S^2{\wedge}S^2{\wedge}S^2 = {\Sigma^3}{\complex}P^2{\wedge}S^2{\wedge}S^2{\wedge}S^2 \cup_{2\nu_{11}} e^{15}$, since $\omega_{n}=2\nu_{n}$ for $n \geq 5$.
An easy computation on the cohomology groups shows that ${\complex}P^2{\wedge}{\complex}P^2$ has the homotopy type of $({\Sigma^2}{\complex}P^2 \vee S^6) \cup_{\beta} e^8$, $\beta$ : $S^7 \overset{\mu}\to S^7{\vee}S^7$ $\xrightarrow{3\nu_{4}\vee\eta}$ $S^4 {\vee} S^6 \subset {\Sigma^2}{\complex}P^2 \vee S^6$, where $\mu$ denotes the unique co-Hopf structure of $S^{7}$.
Then we obtain, up to higher cells in dimension $\geq 10$, that $[({\Sigma^2}{\complex}P^2 \vee S^6) \cup_{\beta} e^8]{\wedge}{\complex}P^2 = ({\Sigma^2}{\complex}P^2{\wedge}{\complex}P^2 \vee {\Sigma}^{6}{\complex}P^2) \cup_{{\Sigma}^{2}\beta} e^{10} = ({\Sigma^2}{\complex}P^2{\wedge}{\complex}P^2 \cup_{3\nu_{6}} e^{10}) \vee {\Sigma}^{6}{\complex}P^2= ({\Sigma^4}{\complex}P^2\cup_{3\nu_{6}}e^{10} \vee S^8) \cup_{{\Sigma}^{2}\beta} e^{10} \vee {\Sigma}^{6}{\complex}P^2 = {\Sigma^4}{\complex}P^2\cup_{3\nu_{6}}e^{10} \vee  {\Sigma^6}{\complex}P^2 \vee {\Sigma}^{6}{\complex}P^2$.
Hence we have, up to higher cells in dimension $\geq 12$, that $[({\Sigma^2}{\complex}P^2 \vee S^6) \cup_{\beta} e^8]({\wedge}{\Sigma^2}{\complex}P^2 \vee S^6) = ({\Sigma^6}{\complex}P^2\cup_{3\nu_{8}}e^{12}) \vee  {\Sigma^8}{\complex}P^2 \vee {\Sigma}^{8}{\complex}P^2 \vee {\Sigma}^{8}{\complex}P^2$.
Thus we obtain that $E^{4}(\G{F_3(\spin{7})})={\Sigma^3}{\complex}P^3{\wedge}S^2{\wedge}S^2{\wedge}S^2 \cup {\Sigma}^{3}{\complex}P^2{\wedge}{\complex}P^2{\wedge}{\complex}P^2{\wedge}{\complex}P^2$ has the homotopy type of 
\begin{align*}&
{\Sigma^3}{\complex}P^3{\wedge}S^2{\wedge}S^2{\wedge}S^2 \cup {\Sigma^3}[({\Sigma^2}{\complex}P^2 \vee S^6) \cup_{\beta} e^8]{\wedge}({\Sigma^2}{\complex}P^2 \vee S^6) 
\\&\hspace{16mm}
\cup ({\Sigma^2}{\complex}P^2 \vee S^6){\wedge}[({\Sigma^2}{\complex}P^2 \vee S^6) \cup_{\beta} e^8] \cup \text{(higher cells $\geq$ $15$)},
\\&
= ({\Sigma}^{9}{\complex}P^2 \cup_{3\nu_{11}} e^{15} \cup_{2\nu_{11}} e^{15}) \vee {\Sigma}^{11}{\complex}P^2 \vee {\Sigma}^{11}{\complex}P^2 \vee {\Sigma}^{11}{\complex}P^2 \cup \text{(higher cells $\geq$ $15$)}
\\&
= ({\Sigma}^{9}{\complex}P^2 \cup_{\nu_{11}} e^{15}) \vee {\Sigma}^{11}{\complex}P^2 \vee {\Sigma}^{11}{\complex}P^2 \vee {\Sigma}^{11}{\complex}P^2 \cup \text{(higher cells $\geq$ $15$)}.
\end{align*}
Then an elementary computation shows that $\pi_{14}(E^{4}(\G{F_3(\spin{7})})) = 0$, and hence $H^{\sigma_3}_3(\alpha) = 0$.
This completes the proof of Theorem \ref{thm:main2}.
\section{Acknowledgements}
The authors thank the University of Aberdeen for its hospitality during their stay in Aberdeen.


\begin{thebibliography}{99}
%
%\bibitem{Borel:hom-cohom-lie}
%A.~Borel, {\em Sur l'homologie et la cohomologie des groupes de Lie compacts connexes}, Amer. J. Math. {\bf 76} (1954), 273--342.
%
%\bibitem{Bott:loop-group}
%R.~Bott, {\em The space of loops on a Lie group}, Michigan Math. J. {\bf 4} (1958), 33--61.
%
\bibitem{FH:category_weight}
E.~Fadell and S.~Husseini, {\em Category weight and Steenrod operations} (Spanish), Bol. Soc. Mat. Mexicana {\bf 37} (1992), 151--161.
%
\bibitem{Ganea:category}
T.~Ganea, {\em Lusternik-Schnirelmann category and strong category}, Illinois. J. Math. {\bf 11} (1967), 417--427.
%
\bibitem{HK:comm_spin}
H.~Hamanaka and A.~Kono, {\em Homotopy-commutativity in spinor groups}, J. Math. Kyoto Univ. {\bf 40} (2000), 389--405.
%
\bibitem{Hardie:fib-cat}
K.~A.~Hardie, {\em A note on fibraitons and category}, Michigan Math. J. \textbf{17}, (1970), 351--352.
%
%\bibitem{Ishitoya:square_spin}
%K.~Ishitoya, {\em Squaring operation in Hermitian symmetric spaces}, J. Math. Kyoto Univ., {\bf 32} (1992), 235--244.
%
\bibitem{IKT:cohomology_spin}
K.~Ishitoya, A.~Kono and H.~Toda, {\em Hopf algebra structure of mod $2$ cohomology of simple Lie groups}, Publ. Res. Inst. Math. Sci. {\bf 12} (1976), 141--167.
%
\bibitem{Iwase:k-ring}
N.~Iwase, {\em On the $K$-ring structure of $X$-projective $n$-space}, Mem. Fac. Sci. Kyushu U. (A) Math., {\bf 38} (1984), 285-297.
%
\bibitem{Iwase:counter-ls}
N.~Iwase, {\em Ganea's conjecture on Lusternik-Schnirelmann category}, Bull. Lon. Math. Soc., {\bf 30} (1998), 623--634.
%
\bibitem{Iwase:counter-ls-m}
N.~Iwase, {\em $A_\infty$-method in Lusternik-Schnirelmann category}, Topology {\bf 41} (2002), 695--723.
%
\bibitem{Iwase:ls-cat-bundle}
N.~Iwase, {\em Lusternik-Schnirelmann category of a sphere-bundle over a sphere}, Topology {\bf 42} (2003), 701--713.
%
\bibitem{Iwase:ls-cat-survay}
N.~Iwase, {\em The Ganea conjecture and recent developments on the Lusternik-Schnirelmann category} (Japanese), S\=ugaku {\bf 56}  (2004), 281--296.
%
\bibitem{IM:higher-associativity}
N.~Iwase and M.~Mimura, {\em Higher homotopy associativity}, Algebraic Topology (Arcata, 1986), 193--220, Lecture Notes in Math. 1370, Springer, Berlin, 1989.
%
\bibitem{IM:cat-sp(3)}
N.~Iwase, M.~Mimura, {L-S categories of simply-connected compact simple Lie groups of low rank}, Lusternik-Schnirelmann category and related topics (South Hadley, MA, 2001), 181--191, Contemp. Math. {\bf 316}, Amer. Math. Soc., Providence, 2002.
%
\bibitem{IMN:spin(7)}
N.~Iwase, M.~Mimura, T.~Nishimoto, {On the cellular decomposition and the Lusternik-Schnirelmann category of $\mathrm{Spin}(7)$}, Top. appl., {\bf 133} (2003), 1--14.
%
\bibitem{IMN:cat-lie}
N.~Iwase, M.~Mimura, T.~Nishimoto, {L-S categories of non-simply-connected compact simple Lie groups}, Topology Appl. {\bf 150}  (2005), 111--123.
%
\bibitem{KK:spinor}
A.~Kono nd K.~Kozima, {\em  The adjoint action of a Lie group on the space of loops},  Japan. J. Math. {\it 45} (1993), 495--510.
%
\bibitem{Rudyak:ls-cat_mfds2}
Y.~B.~Rudyak, {\em On category weight and its applications}, Topology {\bf 38} (1999), 37--55.
%
\bibitem{Stasheff:higher-associativity}
J.~D.~Stasheff, {\em Homotopy associativity of H-spaces, I, II}, Trans. Amer. Math. Soc. {\bf 108} (1963), 275--292, 293--312.
%
\bibitem{Strom:essential-cat-wgt}
J.~Strom, {\em Essential category weight and phantom maps}, Cohomological methods in homotopy theory (Bellaterra, 1998), 409--415, Progr. Math., 196, Birkhauser, Basel, 2001.
%
\bibitem{Varadarajan:fib-cat}
K.~Varadarajan, {\it On fibrations and category}, Math. Z. \textbf{88} (1965), 267--273.
%
\bibitem{Whitehead:elements}
G.~W.~Whitehead, ``Elements of Homotopy Theory'', Springer Verlag, Berlin, GTM series {\bf 61}, 1978.
%
\end{thebibliography}
\end{document}